\documentclass[11pt]{amsart}
\usepackage[letterpaper, margin=1in]{geometry}
\usepackage{amsmath,amssymb,amsthm}
\usepackage{parskip}
\usepackage{tabu}
\usepackage{enumitem}
\usepackage[colorlinks]{hyperref}
\usepackage{thmtools}

\DeclareMathOperator{\St}{St}
\DeclareMathOperator{\Aut}{Aut}
\DeclareMathOperator{\Out}{Out}
\DeclareMathOperator{\Irr}{Irr}
\DeclareMathOperator{\Res}{Res}

\providecommand{\C}{\mathbf C}
\providecommand{\E}{\mathbf E}

\providecommand{\OO}{\mathbf O}
\providecommand{\Z}{\mathbf Z}

\DeclareMathOperator{\M}{M}
\DeclareMathOperator{\J}{J}
\DeclareMathOperator{\Co}{Co}
\DeclareMathOperator{\Fi}{Fi}

\DeclareMathOperator{\PSL}{PSL}
\DeclareMathOperator{\PSU}{PSU}
\DeclareMathOperator{\PSp}{PSp}
\providecommand{\PO}{\mathrm P\Omega}

\declaretheoremstyle[
    spaceabove=-1pt,
    spacebelow=6pt,
    headfont=\normalfont\itshape,
    postheadspace=0.7em,
    qed=\qedsymbol
]{mystyle}
\declaretheorem[name={Proof},style=mystyle,unnumbered]{prf}
\renewenvironment{proof}{\begin{prf}}{\end{prf}}

\newtheorem{theorem}{Theorem}[section]
\newtheorem{proposition}[theorem]{Proposition}
\newtheorem{lemma}[theorem]{Lemma}
\newtheorem{corollary}[theorem]{Corollary}

\hypersetup{urlcolor=blue}

\begin{document}
\title{On the largest Character Degree and Solvable Subgroups of Finite Groups}

\author[Zongshu Wu]{Zongshu Wu}
\address{Stanford Online High School, Redwood City, CA 94063, USA}
\email{zswu07@gmail.com}

\author[Yong Yang]{Yong Yang}
\address{Department of Mathematics, Texas State University, San Marcos, TX 78666, USA}
\email{yang@txstate.edu}

\date{\today}

\begin{abstract}
    Let $G$ be a finite group, and $\pi$ be a set of primes. The $\pi$-core $\OO_\pi(G)$ is the unique maximal normal $\pi$-subgroup of $G$, and $b(G)$ is the largest irreducible character degree of $G$. In 2017, Qian and Yang proved that if $H$ is a solvable $\pi$-subgroup of $G$, then $|H\OO_\pi(G)/\OO_\pi(G)|\le b(G)^3$. In this paper, we improve the exponent of $3$ to $3\log_{504}(168)<2.471$.
\end{abstract}

\maketitle

\section{Introduction}

Let $G$ be a nonabelian finite group, and let $b(G)$ denote the largest irreducible character degree of $G$. We study the relationship between $b(G)$ and the orders of certain solvable subgroups of $G$.

It was proven in \cite{QianShi} that if $P$ is a Sylow $p$-subgroup of $G$, then $|P/\OO_p(G)|<b(G)^2$. This result was generalized in \cite[Theorem 1.2]{QianYang} to arbitrary solvable $\pi$-subgroups, where $\pi$ is a set of primes:

\begin{theorem}
    Let $H$ be a solvable $\pi$-subgroup of $G$. Then $|H\OO_\pi(G)/\OO_\pi(G)|<b(G)^3$.
\end{theorem}

In this paper, we improve the exponent of 3 to $3\log_{504}(168)\approx 2.4703421$, which we will denote by $\mu$ from now on. That is:

\begin{theorem}
    Let $H$ be a solvable $\pi$-subgroup of $G$. Then $|H\OO_\pi(G)/\OO_\pi(G)|\le b(G)^\mu$.
\end{theorem}

The proof will be structured similarly to the original proof of \cite[Theorem 1.2]{QianYang}. To obtain the improved bound, we make use of results in permutation group theory and character theory.

We conjecture that the optimal bound is $b(G)^2$, although we cannot prove it yet.

Throughout this paper, all groups will be assumed to be finite.

\section{Bounds on Finite Simple Groups}

Let $G$ be a group. The automorphism group $\Aut(G)$ naturally acts on $\Irr(G)$, and we shall say that two characters $\chi,\psi\in\Irr(G)$ are \textit{Aut-equivalent} if they belong to the same orbit of this action. In particular, if two characters are Aut-equivalent, then the two corresponding rows in the character table are permutations of each other.

For a group $G$, let $S(G)$ denote the largest order of a solvable subgroup of $G$. We also let $\lambda=\sqrt[3]{24}\approx 2.884499$. Our goal in this section is to prove the following:

\begin{proposition}
    Suppose that $G$ is a nonabelian finite simple group. Then:
    \begin{enumerate}[label = \textup{(\alph*)}]
        \item If $G=\PSL_2(q)$ for a prime power $q\ge 4$, then there exist three non-Aut-equivalent characters $\chi_1,\chi_2,\chi_3\in\Irr(G)$ such that $S(\Aut(G))\le(\chi_1(1)\chi_2(1)\chi_3(1))^{\mu/3}$.
        \item Otherwise, $\lambda\cdot S(\Aut(G))<b(G)^\mu$.
    \end{enumerate}
\end{proposition}
\begin{proof}
    Let $M(G)$ denote the largest order of a proper subgroup of a group $G$, then we have the following inequality:
    \[S(\Aut(G))\le S(G)\cdot|\!\Out(G)|\le M(G)\cdot|\!\Out(G)|.\]
    If $G$ is of Lie type, we bound it from below using the Steinberg character: $b(G)\ge\St(1)$. Note that the value of $\St(1)$ may depend on how $G$ is written: for instance, $\PSU_4(2)\cong\PO_5(3)$, and $\St(1)=64$ for $\PSU_4(2)$ but $\St(1)=81$ for $\PO_5(3)$.

    We first prove part (a): let $G=\PSL_2(q)$ for $q\ge 4$.  We start with the case $q=16$. Using GAP, we find the following three non-Aut-equivalent characters of $\PSL_2(16)$. The following table shows the three rows of the character table, where columns represent conjugacy classes, and $\zeta=e^{2\pi i/15}$.

    \begin{center}\begin{tabular}{c|c|c|c|c|c|c|c|c|c|c}
        & 1 & & & & & & & & \\\hline\hline
        $\chi_1$ & 17 & 1 & $-1$ & 2 & 2 & $-1$ & $-1$ & $-1$ & $-1$ & (8 columns of 0) \\\hline
        $\chi_2$ & 17 & 1 & 2 & $\frac{-1+\sqrt5}2$ & $\frac{-1-\sqrt5}2$ & $\frac{-1+\sqrt5}2$ & $\frac{-1+\sqrt5}2$ & $\frac{-1-\sqrt5}2$ & $\frac{-1-\sqrt5}2$ & (8 columns of 0) \\\hline
        $\chi_3$ & 17 & 1 & $-1$ & $\frac{-1-\sqrt5}2$ & $\frac{-1+\sqrt5}2$ & $\zeta+\zeta^{14}$ & $\zeta^4+\zeta^{11}$ & $\zeta^2+\zeta^{13}$ & $\zeta^7+\zeta^8$ & (8 columns of 0)
    \end{tabular}\end{center}

    These characters all have degree 17, so $(\chi_1(1)\chi_2(1)\chi_3(1))^{\mu/3}>1095$. Using GAP, we also compute that $S(\Aut(G))=960$, so we have $S(\Aut(G))\le(\chi_1(1)\chi_2(1)\chi_3(1))^{\mu/3}$.
    
    Now, we suppose that $q\ne 16$. It is well-known that there exist irreducible characters of degrees $q+1$, $q$, and $q-1$, and we choose these characters, so that $\chi_1(1)\chi_2(1)\chi_3(1)=q^3-q$. If $q$ is one of $5,7,9,11$, then we compute $S(\Aut(G))$ using GAP:

    \begin{center}\begin{tabular}{c|c|c}
        $q$ & $S(\Aut(G))$ & $(q^3-q)^{\mu/3}$ \\\hline\hline
        5 & 24 & $>51$ \\\hline
        7 & 42 & $>120$ \\\hline
        9 & 144 & $>225$ \\\hline
        11 & 110 & $>371$ 
    \end{tabular}\end{center}

    Otherwise, if $q\notin\{5,7,9,11,16\}$, then by \cite{LSU} we have $M(G)=q(q-1)/\gcd(2,q-1)$, therefore $S(\Aut(G))\le fq(q-1)$. We can verify that $fq(q-1)\le(q^3-q)^{\mu/3}$. (In particular, when $q=8$, the bound is tight, as $fq(q-1)=168$ and $q^3-q=504$.)

    Thus, for $G=\PSL_2(q)$, we can choose non-Aut-equivalent characters $\chi_1,\chi_2,\chi_3$ such that
    \[S(\Aut(G))\le(\chi_1(1)\chi_2(1)\chi_3(1))^{\mu/3}.\]

    \bigskip
    
    Now we prove part (b). If $G$ is a sporadic group, then by Chapter 6 of the GAP Character Table Library manual (CTblLibXpls), we know the values of $S(\Aut(G))$. The values of $b(G)$ are well-known, and we obtain the following table:
    \begin{center}\begin{tabular}{c|c|c}
        $G$ & $S(\Aut(G))$ & $b(G)$ \\\hline\hline
        $\M_{11}$ & 144 & 55 \\\hline
        $\M_{12}$ & 432 & 160 \\\hline
        $\M_{22}$ & 1152 & 385 \\\hline
        $\M_{23}$ & 1152 & 2024 \\\hline
        $\M_{24}$ & 13824 & 10395 \\\hline
        $\J_1$ & 168 & 209 \\\hline
        $\J_2$ & 2304 & 336 \\\hline
        $\J_3$ & 3888 & 3078 \\\hline
        $\J_4$ & 28311552 & 3054840657 \\\hline
        HS & 4000 & 3200 \\\hline
        McL & 23328 & 10395 \\\hline
        He & 18432 & 23324 \\\hline
        Ru & 49152 & 118784 \\\hline
        Suz & 279936 & 248832 \\\hline
        O'N & 51840 & 234080 \\\hline
        $\Co_3$ & 69984 & 255024 \\\hline
        $\Co_2$ & 2359296 & 2005875 \\\hline
        $\Co_1$ & 84934656 & 551675124 \\\hline
        HN & 4000000 & 5878125 \\\hline
        Ly & 900000 & 71008476 \\\hline
        Th & 944784 & 190373976 \\\hline
        $\Fi_{22}$ & 10077696 & 2729376 \\\hline
        $\Fi_{23}$ & 3265173504 & 559458900 \\\hline
        $\Fi_{24}'$ & 58773123072 & 160313753600 \\\hline
        B & 29686813949952 & 29823129106907100 \\\hline
        M & 2849934139195390 & 258823477531055064045234375
    \end{tabular}\end{center}
    In all of these cases, $\lambda\cdot S(\Aut(G))<b(G)^\mu$.

    If $G$ is an exceptional group, then by \cite{G2F4}, \cite{E678}, and \cite{Twist}, we obtain the following table. From now on, $q=p^f$ denotes the order of the finite field.
    \begin{center}\begin{tabular}{c|c|c|c}
        $G$ & $M(G)$ & $|\!\Out(G)|$ & $\St(1)$ \\\hline\hline
        $E_6(p^f)$ & $<q^{62}/\gcd(3,q-1)$ & $2f\gcd(3,q-1)$ & $q^{36}$ \\\hline
        $E_7(p^f)$ & $<q^{106}/\gcd(2,q-1)$ & $f\gcd(2,q-1)$ & $q^{63}$ \\\hline
        $E_8(p^f)$ & $<q^{191}$ & $f$ & $q^{120}$ \\\hline
        $F_4(p^f)$ & $<q^{37}$ & $\le 2f$ & $q^{24}$ \\\hline
        $G_2(4)$ & 604800 & 2 & 4096 \\\hline
        $G_2(p^f)$, $q\ne2,4$ & $<q^9$ & $\le 2f$ & $q^6$ \\\hline
        $^3D_4((p^f)^3)$ & $<q^{19}$ & $3f$ & $q^{12}$ \\\hline
        $^2E_6((p^f)^2)$ & $<q^{57}/\gcd(3,q+1)$ & $2f\gcd(3,q+1)$ & $q^{36}$ \\\hline
        $\mathrm{Sz}(2^f)$, $f\ge3$ odd & $<q^3$ & $f$ & $q^2$ \\\hline
        $^2F_4(2)'$ & 11232 & 2 & 2048 \\\hline
        $^2F_4(2^f)$, $f\ge3$ odd & $<q^{16}$ & $f$ & $q^{12}$ \\\hline
        $^2G_2(3^f)$, $f\ge3$ odd & $<q^4$ & $f$ & $q^3$
    \end{tabular}\end{center}
    In all of these cases, $\lambda\cdot M(G)\cdot|\!\Out(G)|<\St(1)^\mu$, hence $\lambda\cdot S(\Aut(G))<b(G)^\mu$.

    If $G$ is a projective symplectic, projective special unitary, or orthogonal group, then by \cite{LSU} and \cite{PO}, we obtain the following table.

    \begin{center}\begin{tabular}{c|c|c|c}
        $G$ & $M(G)$ & $|\!\Out(G)|$ & $\St(1)$ \\\hline\hline
        $\PSp_{2n}(p^f)$ & $<q^{2n^2-n+1}/\gcd(2,q-1)$ & $f\gcd(2,q-1)$ & $q^{n^2}$ \\\hline
        $\PSU_3(p^f)$, $q>2$ & $<q^5/\gcd(3,q+1)$ & $2f\gcd(3,q+1)$ & $q^3$ \\\hline
        $\PSU_4(p^f)$ & $<q^{11}/\gcd(4,q+1)$ & $2f\gcd(4,q+1)$ & $q^6$ \\\hline
        $\PSU_n(p^f)$, $n\ge5$ & $<q^{n^2-2n+2}/\gcd(n,q+1)$ & $2f\gcd(n,q+1)$ & $q^{(n^2-n)/2}$ \\\hline
        $\PO_{2n+1}(p^f)$, $p$ odd & $<q^{2n^2-n+1}/\gcd(2,q-1)$ & $f\gcd(2,q-1)$ & $q^{n^2}$ \\\hline
        $\PO_{2n}^+(p^f)$ & $<q^{2n^2-3n+2}/\gcd(4,q^n+1)$ & $\le 24f$ & $q^{n^2-n}$ \\\hline
        $\PO_{2n}^-(p^f)$ & $<q^{2n^2-3n+2}/\gcd(4,q^n+1)$ & $2f\gcd(4,q^n+1)$ & $q^{n^2-n}$
    \end{tabular}\end{center}
    In all of these cases, $\lambda\cdot M(G)\cdot|\!\Out(G)|<\St(1)^\mu$, hence $\lambda\cdot S(\Aut(G))<b(G)^\mu$.

    If $G=A_n$ for $n\ge7$ (note that $A_5\cong\PSL_2(4)$ and $A_6\cong\PSL_2(9)$), then it is well-known that $S(\Aut(A_n))=S(S_n)\le\lambda^{n-1}$, so it suffices to show that $b(A_n)\ge\lambda^{n/\mu}$.
    
    Let $r=\lfloor n/2\rfloor$, then $S_n$ has an irreducible character of degree $\binom{n/2}r$, which implies that $b(A_n)\ge\frac12\binom{n/2}r$. This is greater than $\lambda^{n/\mu}$ for $n\ge11$. For $n=7,8,9,10$, we have $b(A_n)=35,70,216,567$, respectively, and $b(A_n)\ge\lambda^{n/\mu}$ holds in each case. Therefore, $\lambda\cdot S(\Aut(G))<b(G)^\mu$.

    If $G=\PSL_n(p^f)$ for $n\ge 3$, where $(n,q)\ne(4,2)$ (since $A_8\cong\PSL_4(2)$), then by \cite{LSU},
    \[M(G)=\frac{q^{(n^2-n)/2}}{\gcd(n,q-1)}\cdot\prod_{i=1}^{n-1}(q^i-1).\]
    Since $|\!\Out(G)|=2f\gcd(n,q-1)$ and $\St(1)=q^{(n^2-n)/2}$, we calculate that $\lambda\cdot M(G)\cdot|\!\Out(G)|<\St(1)^\mu$ with the exception of $(n,q)=(3,4)$. If $G=\PSL_3(4)$, we compute using GAP and the GAP Character Table Library that $S(\Aut(\PSL_3(4)))=2304$ and $b(\PSL_3(4))=64$. In any case, $\lambda\cdot S(\Aut(G))<b(G)^\mu$, and the proof of (b) is complete.
\end{proof}

\section{Permutation Groups and Colorings}

In this section, we prove a result on permutation groups.

Suppose that a group $G$ acts faithfully on a finite set $X$. We want to color the elements of $X$ with $3$ distinct colors, such that the number of \textit{color-preserving} permutations (that is, $g\in G$ such that $^gx$ and $x$ have the same color for all $x\in X$) is as small as possible. We will use such a coloring later to assign non-Aut-equivalent characters for a collection of groups.

If $G$ is a $2$-group, then we show below that for some $3$-coloring, the only color-preserving permutation is the identity. More formally:

\begin{theorem}
    Let $G$ be a $2$-group acting faithfully on a finite set $X$. Then there exists a function $f:X\to\{1,2,3\}$ such that there is no nontrivial $g\in G$ with $f(^gx)=f(x)$ for all $x\in X$.
\end{theorem}
\begin{proof}
    We may assume that $G$ is transitive, since a permutation group can be decomposed into transitive groups each acting on an orbit. Then $|X|$ is a power of $2$: let $|X|=2^n$.

    Let us first introduce some terminology. For $n\ge0$ and $i\in\{1,2,3\}$, we define $D_n(i)$ recursively:
    \[D_0(i)=i,\quad D_{n+1}(1)=D_n(1)D_n(2),\quad D_{n+1}(2)=D_n(1)D_n(3),\quad D_{n+1}(3)=D_n(2)D_n(3).\] 
    It easily follows that $D_n(1)<D_n(2)<D_n(3)$ for all $n$. For a function $f$, we also define $P(f)$ as the product of $f(x)$ over all $x$ in the domain of $f$.
    
    We induct on $n$ (where $|X|=2^n$), and prove the following: there exists \textit{three} colorings $f_1,f_2,f_3$ such that for each $i\in\{1,2,3\}$, no nontrivial $g\in G$ is color-preserving, and $P(f_i)=D_n(i)$.
    
    If $n=0$, then $|X|=|G|=1$, and we clearly have the three desired colorings. So assume $n\ge1$.
    
    Let $\Delta_1$ be a maximal block, and $Y$ be a system of blocks containing $\Delta_1$, so that $G$ acts on $Y$. Let $H$ be the kernel of this action, so that $G/H$ acts primitively on $Y$. Since $G/H$ is nilpotent, it is well-known that we must have $|G/H|=|Y|=2$. Let $\Delta_2=X\setminus\Delta_1$, so that $\Delta_2$ is also a block.

    For $j\in\{1,2\}$, let $K_j$ be the subgroup of $g\in G$ that fix all elements of $\Delta_j$. We can easily see that $K_j\lhd H$, and that $H/K_j$ acts faithfully and transitively on $\Delta_j$. Therefore, by the inductive hypothesis, we may find three colorings $f_{1j},f_{2j},f_{3j}:\Delta_j\to\{1,2,3\}$ such that for each coloring, the set of color-preserving $g\in H$ is $K_j$. That is, any color-preserving $g\in H$ acts trivially on $\Delta_j$.
    
    We combine these colorings as follows:
    \[f_1(x)=\left\{\negthickspace\!\begin{array}{ll}
        f_{11}(x)&x\in\Delta_1 \\ f_{22}(x)&x\in\Delta_2
    \end{array}\right.,\quad f_2(x)=\left\{\negthickspace\!\begin{array}{ll}
        f_{11}(x)&x\in\Delta_1 \\ f_{32}(x)&x\in\Delta_2
    \end{array}\right.,\quad f_3(x)=\left\{\negthickspace\!\begin{array}{ll}
        f_{21}(x)&x\in\Delta_1 \\ f_{32}(x)&x\in\Delta_2
    \end{array}\right..\]
    Let $i\in\{1,2,3\}$. We may easily verify that $P(f_i)=D_n(i)$. Now suppose that some permutation $g\in G$ is color-preserving. If $g$ is not an element of $H$, then
    \[\prod_{x\in\Delta_1}f_i(x)
    =\prod_{x\in\Delta_1}f_i(^gx)
    =\prod_{x\in\Delta_1}f_i(^gx)
    =\prod_{x\in\Delta_2}f_i(x),\]
    a contradiction, since $P(f_{1j})\ne P(f_{2k})$ for $j\ne k$. Thus $g\in H$, and since it preserves the coloring of $\Delta_1$ and $\Delta_2$, we know that $g$ acts trivially on $\Delta_1$ and $\Delta_2$, hence $g$ is the identity.
\end{proof}

This result allows us to generalize to all solvable groups: if $G$ is solvable, then for some $3$-coloring, the number of color-preserving permutations of $G$ is most $\sqrt{|G|}$.

\begin{corollary}
    Let $G$ be a solvable group acting faithfully on a finite set $X$. Then there exists a function $f:X\to\{1,2,3\}$ such that $|J_f|^2\le|G|$, where we define $J_f=\{g\in G:\forall x\in X,f(^gx)=f(x)\}$. Furthermore, we may have $|f^{-1}(1)|\ge|f^{-1}(2)|\ge|f^{-1}(3)|$.
\end{corollary}
\begin{proof}
    Let $G_1$ be a Sylow $2$-subgroup of $G$, and since $G$ is solvable, we can let $G_2$ be a Hall $2'$-subgroup of $G$. We easily see that $|G|=|G_1||G_2|$.

    If $|G_1|\ge|G_2|$, then we apply Theorem 3.1 on $G_1$ to obtain a coloring $f:X\to\{1,2,3\}$ such that $J_f\cap G_1$ is trivial. It follows that
    \[|J_f|=|J_fG_1|/|G_1|\le|G|/|G_1|\le\sqrt{|G|}.\]
    If $|G_2|\ge|G_1|$, then since $|G_2|$ is odd, Gluck's permutation lemma (see \cite[Corollary 5.7]{ManzWolf}) implies the existence of a set $S\subseteq X$ such that no nontrivial element of $G_2$ fixes $S$. We set $f(x)=1$ if $x\in S$, and $f(x)=2$ otherwise, so that $J_f\cap G_2$ is trivial. Similarly to above, we obtain $|J_f|\le\sqrt{|G|}$.

    We may ensure that $|f^{-1}(1)|\ge|f^{-1}(2)|\ge|f^{-1}(3)|$ by simply permuting $1,2,3$.
\end{proof}

\section{Proof of the Main Result}

In this section, we finish the proof of Theorem 1.2. The structure will be similar to the proof of \cite[Theorem 1.2]{QianYang}, but we will use a slightly different framework.

First, let us briefly recall how to decompose a subgroup of a product group.

\begin{lemma}
    Let $G=X\rtimes Y$ and $H\le G$. Then $H/(H\cap X)$ is isomorphic to a subgroup of $Y$.
\end{lemma}
\begin{proof}
    By the Second Isomorphism Theorem, $H/(H\cap X)\cong HX/X\le G/X\cong Y$.
\end{proof}

This lemma allows us to decompose $H$ into $H\cap X$ and $H/(H\cap X)$. In particular, suppose that $H$ is a subgroup of $G=N_1\times N_2\times\dots\times N_k$. Then by repeatedly applying Lemma 4.1, we may decompose $H$ into $H_1,H_2,\dots,H_k$, where $H_i\le N_i$ and $|H|=|H_1||H_2|\cdots|H_k|$.

An easy consequence of this is that:

\begin{lemma}
    Let $G=N_1\times N_2\times\dots\times N_k$, and let $H$ be a solvable subgroup of $G$. Then
    \[|H|\le\prod_{i=1}^kS(N_i).\]
\end{lemma}

\bigskip

Finally, we begin the proof of our main result.

\textit{Proof of Theorem 1.2.} Since $b(G/\OO_\pi(G))\le b(G)$, we may assume that $\OO_\pi(G)=1$.

Recall that $\E(G)$ is the layer of $G$. Let $Z=\Z(\E(G))$, and $K=G/Z$. If any $p\in\pi$ divides $|Z|$, then the Sylow $p$-subgroup of $Z$ is characteristic in $G$, contradicting $\OO_\pi(G)=1$. Thus, $Z$ is a $\pi'$-group, and hence intersects trivially with $H$.

Observe that $H$ is isomorphic to $HZ/Z\le K$, which we will denote by $\bar H$.

\bigskip

It is well-known that $\E(K)=\E(G)/Z$, so it decomposes into a direct product of nonabelian simple groups. Let $\E(K)=E_1\times E_2\times\dots\times E_m$, where $E_i\le\E(K)$ is isomorphic to $L_i^{k_i}$ for some nonabelian simple group $L_i$, such that the $L_i$ are pairwise nonisomorphic. We also write $E_i=L_{i1}\times\cdots\times L_{ik_i}$, where $L_{ij}\le E_i$ is isomorphic to $L_i$.

Let $X_i$ denote the set of $L_{ij}$, and let $X$ denote the union of the $X_i$. Let $C=\C_K(\E(K))$ and $A=K/C$, so that $A$ acts faithfully on $\E(K)$ via conjugation. This induces an action of $A$ on each $X_i$, as any $L_{ij}$ must be mapped to $L_{il}$ for some $l$. Of course, we also have an action on $X$ itself.

Let $N$ be the kernel of the action on $X$, which is easily seen to be normal in $A$. Note that $N$ is isomorphic to a subgroup of $\Aut(L_1)^{k_1}\times\dots\times\Aut(L_m)^{k_m}$. Next, let $P=A/N$. The action of $A$ on $X$ induces a faithful action of $P$ on $X$ and (not necessarily faithful) actions on each $X_i$, and therefore $P$ is isomorphic to a subgroup of $S_{k_1}\times\dots\times S_{k_m}$.

\bigskip

Let $H_C=\bar H\cap C\le C$, $H_N=(\bar HC/C)\cap N\le N$, and $H_P=(\bar HC/C)N/N\le P$, so that
\[|H|=|H_C||H_N||H_P|.\]
Using the ``large orbit theorem'' \cite[Corollary 1.2]{Dolfi}, it was shown in the proof of \cite[Theorem 3.2]{QianYang} (where $C$ and $H_C$ were referred to as ``$N_1$'' and ``$H_0$'') that
\[|H_C|\le b(C)^2.\]
Since $H_N$ embeds into $\Aut(L_1)^{k_1}\times\dots\times\Aut(L_m)^{k_m}$, by Lemma 4.2 we have
\[|H_N|\le\prod_{i=1}^mS(\Aut(L_i))^{k_i}.\]
Similarly, since $H_P$ embeds into $S_{k_1}\times\dots\times S_{k_m}$, we may decompose it into $H_1,\dots,H_m$, such that each $H_i$ is solvable and acts faithfully on $X_i$. 

\bigskip

To proceed, we carefully choose characters $\theta_{ij}\in\Irr(L_{ij})$ in the manner described below. Once we have fixed the choices of $\theta_{ij}$, let $\theta_i=\theta_{i1}\times\dots\times\theta_{ik_i}\in\Irr(E_i)$.

If $L_i\cong\PSL_2(q)$ for some $q\ge4$, we apply Corollary 3.2 to the action of $H_i$ on $X_i$ to obtain a function $f_i:X_i\to\{1,2,3\}$, such that if $J_i=\{h\in H_i:\forall j,f_i(^hL_{ij})=f_i(L_{ij})\}$, then $|J_i|^2\le|H_i|$. We also have $|f_i^{-1}(1)|\ge|f_i^{-1}(2)|\ge|f_i^{-1}(3)|$.

Proposition 2.1 (a) provides us three non-Aut-equivalent characters $\chi_1,\chi_2,\chi_3\in\Irr(L_i)$, such that $S(\Aut(L_i))\le(\chi_1(1)\chi_2(1)\chi_3(1))^{\mu/3}$.
We may assume that $\chi_1(1)\ge\chi_2(1)\ge\chi_3(1)$. Choose $\theta_{ij}$ to be $\chi_{f_i(L_{ij})}$. Then
\[\theta_i(1)=\chi_1(1)^{|f_i^{-1}(1)|}\cdot\chi_2(1)^{|f_i^{-1}(2)|}\cdot\chi_3(1)^{|f_i^{-1}(3)|}\ge(\chi_1(1)\chi_2(1)\chi_3(1))^{k_i/3}.\]
Hence,
\[S(\Aut(L_i))^{k_i}\cdot|H_i|\le(\chi_1(1)\chi_2(1)\chi_3(1))^{\mu k_i/3}\cdot[H_i:J_i]^2\le\theta_i(1)^\mu\cdot[H_i:J_i]^\mu.\]
Otherwise, if $L_i$ is not a $\PSL_2(q)$, then set $\theta_{ij}$ to be a character of largest degree, so that $\theta_i$ has degree $b(L_i)^{k_i}$. We let $J_i=H_i$, then by Proposition 2.1 (b), we also have
\[S(\Aut(L_i))^{k_i}\cdot|H_i|\le S(\Aut(L_i))^{k_i}\cdot\lambda^{k_i}<b(L_i)^{\mu k_i}=\theta_i(1)^\mu\cdot[H_i:J_i]^\mu.\]
Since $\E(K)$ is centerless, it intersects trivially with $C$, so we can form the direct product $\E(K)\times C$. Let $\psi\in\Irr(C)$ such that $\psi(1)=b(C)$, and let $\theta=\theta_1\times\dots\times\theta_m\times\psi\in\Irr(\E(K)\times C)$. Then by multiplying together the inequalities above, along with the bound on $|H_C|$, we obtain:
\[|H|=|H_C||H_N||H_P|\le\psi(1)^\mu\cdot\prod_{i=1}^m\theta_i(1)^\mu\cdot\prod_{i=1}^m[H_i:J_i]^\mu=\theta(1)^\mu\cdot\prod_{i=1}^m[H_i:J_i]^\mu.\]

\bigskip

By Frobenius reciprocity, there exists a character $\chi\in\Irr(K)$ such that $\theta\in\Irr(\E(K)\times C)$ is an irreducible constituent of $\Res(\chi)$, the restriction of $\chi$ to the subgroup $\E(K)\times C$. By Clifford's theorem, the degree of this character is
\[\chi(1)=\langle\theta,\Res(\chi)\rangle\cdot[K:I_K(\theta)]\cdot\theta(1),\]
where recall that $I_K(\theta)$ is the inertia group of $\theta$.

Let $I=((I_K(\theta)\cap\bar H)C/C)N/N\le H_P$. We can easily show that $[K:I_K(\theta)]\ge[H_P:I]$. In the same way as with $H_P$, we may decompose $I$ into $I_1,I_2,\dots,I_m$ such that $I_i\le H_i$.

Note that each $I_i$ also acts on $X_i$. We claim that $I_i\le J_i$. If $L_i$ is not a $\PSL_2(q)$, then $J_i=H_i$ and the claim is trivial. Otherwise, let $h\in I_i$. Let us write $^hL_{ij}=L_{i\sigma_h(j)}$. Then, since an isomorphism sends $\theta_{ij}$ to $\theta_{i\sigma_h(j)}$, these characters must be Aut-equivalent, which implies that $f_i(L_{ij})=f_i(L_{i\sigma_h(j)})=f_i(^hL_{ij})$ for all $j$. That is, $h\in J_i$. Hence, $I_i\le J_i$. Therefore, it follows that
\[|H|\le\theta(1)^\mu\cdot\prod_{i=1}^m[H_i:I_i]^\mu=\theta(1)^\mu\cdot[H_P:I]\le\chi(1)^\mu\le b(K)^\mu\le b(G)^\mu,\]
as desired. $\hfill\square$

\section{Further Discussions}

In this section, we discuss potential improvements to the bound $|H\OO_\pi(G)/\OO_\pi(G)|\le b(G)^\mu$.

First of all, Proposition 2.1 (a) is tight: if $G=\PSL_2(8)$, then using GAP, we find that the five characters of largest degrees, shown in the following table, where $\zeta=e^{2\pi i/7}$:

\begin{center}\begin{tabular}{c|c|c|c|c|c|c|c|c|c}
    & 1 & & & & & & & \\\hline\hline
    $\psi_1$ & 7 & $-1$ & $-2$ & 0 & 0 & 0 & 1 & 1 & 1 \\\hline
    $\psi_2$ & 8 & 0 & $-1$ & 1 & 1 & 1 & $-1$ & $-1$ & $-1$ \\\hline
    $\psi_3$ & 9 & 1 & 0 & $\zeta+\zeta^6$ & $\zeta^2+\zeta^5$ & $\zeta^3+\zeta^4$ & 0 & 0 & 0 \\\hline
    $\psi_4$ & 9 & 1 & 0 & $\zeta^2+\zeta^5$ & $\zeta^3+\zeta^4$ & $\zeta+\zeta^6$ & 0 & 0 & 0 \\\hline
    $\psi_5$ & 9 & 1 & 0 & $\zeta^3+\zeta^4$ & $\zeta+\zeta^6$ & $\zeta^2+\zeta^5$ & 0 & 0 & 0
\end{tabular}\end{center}

By using GAP, we can show that the three characters of degree 9 are Aut-equivalent. Therefore, for any non-Aut-equivalent characters $\chi_1,\chi_2,\chi_3\in\Irr(G)$, we have $\chi_1(1)\chi_2(1)\chi_3(1)\le7\cdot8\cdot9=504$. Furthermore, using GAP we also compute that $S(\Aut(G))=168$. So $\mu=3\log_{504}(168)$ is the smallest constant that makes Proposition 2.1 work.

Now we examine if there is room for improvement in the main proof. Suppose that $\OO_\pi(G)=1$ and the equality $|H|=b(G)^\mu$ is attained. Then, using the same notation as in the main proof, we see that the following facts must hold true:
\begin{enumerate}
    \item $H_C$ is trivial and $b(C)=1$;
    \item $\E(K)$ is isomorphic to $\PSL_2(8)^k$ for $k\ge 0$;
    \item If $\E(K)$ is nontrivial, then $I_1=J_1=H_1$.
\end{enumerate}

Unfortunately, we are unable to immediately arrive at a contradiction, but these restrictions on the group are quite strong, and may be useful to obtain improvements in the future.

Finally, we conjecture the following:

\textbf{Conjecture.} Let $H$ be a solvable $\pi$-subgroup of $G$. Then $|H\OO_\pi(G)/\OO_\pi(G)|\le b(G)^2$.

We believe that this is true, because it holds for all finite nonabelian simple groups: in fact, we can show that $S(G)<b(G)^2$, where $G$ is a nonabelian finite simple group. The proof of this is similar to that of Proposition 2.1, and we omit it.

\section{Acknowledgements}
Yang was partially supported by a grant from the Simons Foundation (\#918096, to YY).

\end{document}